\documentclass[12pt,oneside,reqno]{amsart}
\usepackage{mathrsfs}
\usepackage{graphics}
\usepackage{amssymb}
\usepackage{enumerate}
\newcommand{\dif}{\mathrm{d}}

\newcommand{\be}{\begin{eqnarray}}
\newcommand{\ee}{\end{eqnarray}}
\newcommand{\ce}{\begin{eqnarray*}}
\newcommand{\de}{\end{eqnarray*}}
\newtheorem{theorem}{Theorem}[section]
\newtheorem{lemma}[theorem]{Lemma}
\newtheorem{remark}[theorem]{Remark}
\newtheorem{definition}[theorem]{Definition}
\newtheorem{proposition}[theorem]{Proposition}
\newtheorem{Example}[theorem]{Example}
\newtheorem{corollary}[theorem]{Corollary}
\def\e{\varepsilon}
\def\t{\theta}
\def\a{\alpha}

\def\b{\beta}

\def\[{{\Big[}}
\def\]{{\Big]}}
\def\<{{\langle}}
\def\>{{\rangle}}
\def\({{\Big(}}
\def\){{\Big)}}

\def\tr{{\rm tr}}

\def\no{\nonumber}
\def\bt{\begin{theorem}}
\def\et{\end{theorem}}
\def\bl{\begin{lemma}}
\def\el{\end{lemma}}
\def\br{\begin{remark}}
\def\er{\end{remark}}
\def\bx{\begin{Example}}
\def\ex{\end{Example}}
\def\bd{\begin{definition}}
\def\ed{\end{definition}}
\def\bp{\begin{proposition}}
\def\ep{\end{proposition}}
\def\bc{\begin{corollary}}
\def\ec{\end{corollary}}

\def\cB{{\mathcal B}}
\def\cC{{\mathcal C}}

\def\cK{{\mathcal K}}
\def\cL{{\mathcal L}}
\def\cM{{\mathcal M}}

\def\cP{{\mathcal P}}

\def\cT{{\mathcal T}}

\def\mE{{\mathbb E}}

\def\mN{{\mathbb N}}

\def\mP{{\mathbb P}}
\def\mQ{{\mathbb Q}}
\def\mR{{\mathbb R}}
\def\mS{{\mathbb S}}

\def\mW{{\mathbb W}}

\def\sB{{\mathscr B}}
\def\sC{{\mathscr C}}

\def\sF{{\mathscr F}}

\def\sL{{\mathscr L}}

\def\sS{{\mathscr S}}

\def\geq{\geqslant}
\def\leq{\leqslant}

\begin{document}

\allowdisplaybreaks

\title{Superposition principles for the Zakai equations and the Fokker-Planck equations on measure spaces*}

\author{Huijie Qiao}

\dedicatory{School of Mathematics,
Southeast University\\
Nanjing, Jiangsu 211189,  China\\
Department of Mathematics, University of Illinois at
Urbana-Champaign\\
Urbana, IL 61801, USA\\
hjqiaogean@seu.edu.cn}



\thanks{*This work was partly supported by NSF of China (No. 11001051, 11371352) and China Scholarship Council under Grant No. 201906095034.}


\date{}

\begin{abstract}
The work concerns the superposition between the Zakai equations and the Fokker-Planck equations on measure spaces. First, we prove a superposition principle for the Fokker-Planck equations on $\mR^\mN$ under the integrable condition. And then by means of it, we show two superposition principles for the weak solutions of the Zakai equations from the nonlinear filtering problems and the weak solutions of the Fokker-Planck equations on measure spaces. As a by-product, 
we give some weak conditions under which the Fokker-Planck equations can be solved in the weak sense.
\end{abstract}

\maketitle \rm

\section{Introduction}

Nonlinear filtering problems are to extract information about the signal processes from the observation processes and then estimate and predict the signal processes. And nonlinear filtering problems have been widely applied to various fields, such as physics, biology, the control theory and the weather forecast. Moreover, more and more researchers are paying attention to nonlinear filtering problems.

Nonlinear filtering  problems are closely related with two types of measure-valued equations-the Zakai equations and the Kushner-Stratonovich equations. Furthermore, the Zakai equations are linear, and the Kushner-Stratonovich equations are not linear. Hence, in order to solve nonlinear filtering problems, many researchers usually deduce and study the Zakai equations. Let us mention some results related with our work. In \cite{sz}, Szpirglas took the Zakai equations as stochastic differential equations 
and studied the uniqueness of weak solutions for the Zakai equations. There the signal processes are Markov processes independent of the Wiener processes in the observation processes. In \cite{lh} Lucic-Heunis also investigated the uniqueness problem 
of weak solutions for the Zakai equations when the signal processes depend on the Wiener process in the observation processes. Later, the author  proved the uniqueness of weak solutions to the Zakai equations for non-Gaussian signal-observation 
systems with independent noises in \cite{q3} and with correlated noises in \cite{q2}, respectively.

In the paper, we still view the Zakai equations as stochastic differential equations and set up a correspondence between the weak solutions of the Zakai equations and the weak solutions of the Fokker-Planck equations on some measure space. 
Concretely speaking, fix $T>0$ and consider the following signal-observation system $(X_t, Y_t)$ on $\mR^n\times\mR^m$:
\be\left\{\begin{array}{l}
\dif X_t=b_1(t,X_t)\dif t+\sigma_0(t,X_t)\dif B_t+\sigma_1(t,X_t)\dif W_t,\\
\dif Y_t=b_2(t,X_t)\dif t+\sigma_2(t)\dif W_t, \quad 0\leq t\leq T,
 \end{array}
\right. 
\label{Eq1} 
\ee
where $B, W$ are $d$-dimensional and $m$-dimensional Brownian motions defined on a complete filtered probability space $(\Omega, \mathscr{F}, \{\mathscr{F}_t\}_{t\in[0,T]},\mP)$, respectively. Moreover, $W_{\cdot}, B_{\cdot}$ are mutually independent. The mappings $b_1:[0,T]\times\mR^n\mapsto\mR^n$, $\sigma_0:[0,T]\times\mR^n\mapsto\mR^{n\times d}$, $\sigma_1:[0,T]\times\mR^n\mapsto\mR^{n\times m}$ and $b_2:[0,T]\times\mR^n\mapsto\mR^m$, $\sigma_2:[0,T]\mapsto\mR^{m\times m}$ are all Borel measurable. The initial value $X_0$ is assumed to be a square integrable random variable independent of $Y_0, B_{\cdot}, W_{\cdot}$. The nonlinear filtering problem related with the system (\ref{Eq1}) is usually called as the nonlinear filtering problem 
with correlated noises. And then we deduce the Zakai equation (\ref{zakaieq0}) (c.f. Section \ref{Zakaiequ}) associated with the system (\ref{Eq1}) and define its weak solution. After this, we construct the Fokker-Planck equation (\ref{fpecoup}) (c.f. Section \ref{TFPE}) related with (\ref{zakaieq0}) and prove a superposition principle between Eq.(\ref{zakaieq0}) and Eq.(\ref{fpecoup}).

Moreover, we also consider the other signal-observation system $(\check{X}_t, \check{Y}_t)$ on $\mR^n\times\mR^m$:
\be\left\{\begin{array}{l}
\dif \check{X}_t=\check{b}_1(t,\check{X}_t)\dif t+\check{\sigma}_1(t,\check{X}_t)\dif W_t,\\
\dif \check{Y}_t=\check{b}_2(t,\check{X}_t)\dif t+\check{\sigma}_2\dif W_t+\check{\sigma}_3\dif B_t, \quad 0\leq t\leq T,
 \end{array}
\right. 
\label{Eq0} 
\ee
where $W, B$ are the same to ones in the system (\ref{Eq1}). The initial value $\check{X}_0$ is assumed to be a square integrable random variable independent of $\check{Y}_0, W_{\cdot}, B_{\cdot}$. The mappings $\check{b}_1:[0,T]\times\mR^n\mapsto\mR^n$, $\check{\sigma}_1:[0,T]\times\mR^n\mapsto\mR^{n\times m}$ and $\check{b}_2:[0,T]\times\mR^n\mapsto\mR^m$ are all Borel measurable. $\check{\sigma}_2, \check{\sigma}_3$ are $m\times m$ and $m\times d$ real matrices, respectively. The nonlinear filtering problem related with the system (\ref{Eq0}) is usually called as the nonlinear filtering problem with correlated sensor noises. We also give the Zakai equation (\ref{zakaieq01}) (c.f. Section \ref{Zakaiequ}) associated with the system (\ref{Eq0}) and the Fokker-Planck equation (\ref{fpecoup2}) (c.f. Section \ref{TFPE}) related with (\ref{zakaieq01}) and show a superposition principle between Eq.(\ref{zakaieq01}) and Eq.(\ref{fpecoup2}). It is worthwhile to mentioning that if $\sigma_1=0$ in the system (\ref{Eq1}) and $\check{\sigma}_2=0$ in the system (\ref{Eq0}), both the systems (\ref{Eq1}) and (\ref{Eq0}) are just the usual signal-observation systems. That is, the systems (\ref{Eq1}) and (\ref{Eq0}) cover a lot of signal-observation systems.
 
Here is a summary of our results. First, we prove a superposition principle for the Fokker-Planck equations on $\mR^\mN$ under the integrable condition (See Theorem \ref{supeprinrn}). Theorem \ref{supeprinrn} generalizes \cite[Theorem 7.1]{ad} for the continuity equations, \cite[Theorem 7.1]{T0} for the $p$-order integrable condition and \cite[Theorem 2.5]{T} for the finite dimensional space. Second, we show two superposition principles for the weak solutions of the Zakai equations from the nonlinear filtering problems and the weak solutions of the Fokker-Planck equations on some measure space (See Theorem \ref{supeprin} and \ref{supeprin2}). As far as we know, this is the first time to set up a correspondence between the Zakai equations and the Fokker-Planck equations on measure spaces. Moreover, due to the speciality of the Zakai equations, by the way some complicated Fokker-Planck equations can be solved. Therefore, these results are useful for application of Fokker-Planck equations. It is interesting to study the correspondence between the Zakai equations from nonlinear filtering problems of non-Gaussian systems  and the nonlocal Fokker-Planck equations on measure spaces as \cite{q10}. This is our forthcoming work.

In Section \ref{pre}, we introduce some notation and define L-derivatives for functions on $\cP_2(\mR^n)$. And then we prove a superposition principle for Fokker-Planck equations on $\mR^\mN$ in Section \ref{susupeprinin}. In Section \ref{Zakaiequ}, two types of Zakai equations from nonlinear filtering problems are deduced. And then, we give two types of Fokker-Planck equations on $\cM(\mR^n)$. Finally, we prove two superposition principles 
for the Zakai equations and Fokker-Plank equations on $\cM(\mR^n)$ in Section \ref{susupeprin}. 

The following convention will be used throughout the paper: $C$ with or without indices will denote
different positive constants whose values may change from one place to another.

\section{Preliminary}\label{pre}

In the section, we introduce some notation and define L-derivatives for functions on $\cP_2(\mR^n)$.

\subsection{Notation}

In this subsection, we introduce some  notation used in the sequel. 

For convenience, we shall use $\mid\cdot\mid$
and $\parallel\cdot\parallel$  for norms of vectors and matrices, respectively. Let $A^*$ denote the transpose of the matrix $A$.

Let $\sB(\mR^n)$ be the Borel $\sigma$-field on $\mR^n$. Let $\cB_b(\mR^n)$ denote the set of all real-valued uniformly bounded $\mathscr{B}(\mR^n)$-measurable functions on $\mR^n$. $C^2(\mR^n)$ stands for the space of continuous functions on $\mR^n$ which have continuous partial derivatives of order up to $2$, and $C_b^2(\mR^n)$ stands for the subspace of $C^2(\mR^n)$, consisting of functions whose derivatives up to order 2 are bounded. $C_c^2(\mR^n)$ is the collection of all functions in $C^2(\mR^n)$ with compact support and $C_c^\infty(\mR^n)$ denotes the collection of all real-valued $C^\infty$ functions of compact
support.

Let $\cM(\mR^n)$ be the set of bounded Borel measures on $\sB(\mR^n)$. And then we endow $\cM(\mR^n)$ with the weak convergence topology. Let $\cP(\mR^n)$ be the space of all probability measures in $\cM(\mR^n)$. $\cP_2(\mR^n)$ denotes the subspace of $\cP(\mR^n)$ with finite second moments. And the distance of $\mu, \nu\in\cP_2(\mR^n)$ is defined as 
$$
\mW^{2}_{2}(\mu,\nu):=\inf_{\pi\in\sC(\mu_1, \mu_2)} \int_{\mR^n \times \mR^n} |x-y|^2 \pi(\dif x, \dif y),
$$
where $\sC(\mu_1, \mu_2)$ denotes the set of  all the probability measures whose marginal distributions are $\mu_1, \mu_2$, respectively. Thus, $(\cP_2(\mR^n),\mW_{2})$ is a Polish space.

\subsection{L-derivatives for functions on $\cP_2(\mR^n)$}\label{lde} 

In the subsection we recall the definition of L-derivatives for functions on $\cP_2(\mR^n)$. 

The definition of L-derivatives was first introduced by Lions \cite{Lion}. Moreover, he used some abstract probability spaces to describe the L-derivatives. Here, for the convenience to understand the definition, we apply a straight way to state it (\cite{rw}). Let $I$ be the identity mapping on $\mR^n$. For $\mu\in\cP_2(\mR^n)$ and $\varphi\in L^2(\mR^n, \sB(\mR^n), \mu;\mR^n)$, by simple calculation, it holds that $\mu\circ(I+\varphi)^{-1}\in\cP_2(\mR^n)$. $<\mu, \varphi>:=\int_{\mR^n}\varphi(x)\mu(\dif x)$.

\bd\label{lderi}
(i) A function $H: \cP_2(\mR^n)\mapsto\mR$ is called L-differentiable at $\mu\in\cP_2(\mR^n)$, if the functional 
$$
L^2(\mR^n, \sB(\mR^n), \mu;\mR^n)\ni\varphi\mapsto H(\mu\circ(I+\varphi)^{-1})
$$
is Fr\'echet differentiable at $0\in L^2(\mR^n, \sB(\mR^n), \mu;\mR^n)$; that is, there exists a unique $\gamma\in L^2(\mR^n, \sB(\mR^n), \mu;\mR^n)$ such that 
$$
\lim\limits_{<\mu,|\varphi|^2>\rightarrow 0}\frac{H(\mu\circ(I+\varphi)^{-1})-H(\mu)-\mu(\gamma\cdot\varphi)}{\sqrt{<\mu, |\varphi|^2>}}=0.
$$
In the case, we denote $\partial_{\mu}H(\mu):=\gamma$ and call it the L-derivative of $H$ at $\mu$.

(ii) A function $H: \cP_2(\mR^n)\mapsto\mR$ is called L-differentiable on $\cP_2(\mR^n)$ if the L-derivative $\partial_{\mu}H(\mu)$ exists for all $\mu\in\cP_2(\mR^n)$.

(iii) By the same way, $\partial^2_\mu H(\mu)(y,y')$ for $y, y'\in\mR^n$ can be defined.
\ed

Next, we give an example to explain how to compute the L-derivatives for cylindrical functions on $\cP_2(\mR^n)$. For example, $H(\mu):=g(<\mu,\varphi_1>, <\mu,\varphi_2>, \cdots, <\mu,\varphi_k>), k\in\mN, g\in C_b^2(\mR^{k})$, $\mu\in\cP_2(\mR^n)$, $\varphi_1, \cdots, \varphi_k\in C_c^2(\mR^n)$. And then by simple calculation, we know that for $y, y'\in\mR^n$
\be
&&\partial_{\mu}H(\mu)(y)=\partial_ig(<\mu,\varphi_1>, <\mu,\varphi_2>, \cdots, <\mu,\varphi_k>)\partial_y\varphi_i(y), \no\\
&&\partial_y\partial_{\mu}H(\mu)(y)=\partial_ig(<\mu,\varphi_1>, <\mu,\varphi_2>, \cdots, <\mu,\varphi_k>)\partial^2_y\varphi_i(y), \no\\
&&\partial^2_{\mu}H(\mu)(y,y')=\partial_{ij}g(<\mu,\varphi_1>, <\mu,\varphi_2>, \cdots, <\mu,\varphi_k>)\partial_y\varphi_i(y)\partial_{y'}\varphi_j(y').\label{cylilderi}
\ee
Here and hereafter, we use the convention that repeated indices imply summation. 

\section{A superposition principle for the Fokker-Planck equations on $\mR^\mN$}\label{susupeprinin}

In the section, we prove a superposition principle for Fokker-Planck equations (FPEs in short) on $\mR^\mN$.

Firstly, we write $\mR^\mN$ as $\mR^\infty$ and endow $\mR^\infty$ with the product topology. And then for $\boldsymbol{x}, \boldsymbol{y}\in\mR^\infty$, define
$$
d_{\mR^\infty}(\boldsymbol{x}, \boldsymbol{y}):=\sum\limits_{k=1}^\infty 2^{-k}\min\{1,|\boldsymbol{x}^k-\boldsymbol{y}^k|\}.
$$
So, $(\mR^\infty,d_{\mR^\infty})$ is a complete and separable metric space. Accordingly, we consider the space $C([0,T],\mR^\infty)$ endowed with the distance
$$
d_{C([0,T],\mR^\infty)}(w, \tilde{w}):=\sum\limits_{k=1}^\infty 2^{-k}\max\limits_{t\in[0,T]}\min\{1,|w_t^k-\tilde{w}_t^k|\}.
$$
Thus, $(C([0,T],\mR^\infty),d_{C([0,T],\mR^\infty)})$ also becomes a complete and separable metric space. A compact set $\cK\subset C([0,T],\mR^\infty)$ means that the set $\{\boldsymbol{x}^k\circ w: w\in\cK\}$ is compact in $C([0,T],\mR)$ for $k\in\mN$. 
Hence, if $\Psi^k: C([0,T],\mR)\rightarrow [0,\infty]$ is a coercive functional for $k\in\mN$ (i.e. $\{\Psi^k\leq C\}$ is compact for every $C\geq 0$), we say that $\Psi: C([0,T],\mR^\infty)\rightarrow [0,\infty]$ defined by
$$
\Psi(w):=\sum\limits_{k=1}^\infty\Psi^k(\boldsymbol{x}^k\circ w), \qquad w\in C([0,T],\mR^\infty)
$$
is coercive. Let $C_T^\infty:=C([0,T],\mR^\infty)$. And $w$ stands for a generic element in $C_T^\infty$. For any $t\in[0,T]$, set 
$$
e_t: C_T^\infty\rightarrow\mR^\infty, \quad e_t(w)=w_t, \quad w\in C_T^\infty.
$$
Let $\cB_t:=\sigma\{w_s: s\in[0,t]\}$,  $\bar{\cB}_t:=\cap_{s>t}\cB_s$, and $\cB:=\cB_T$. 

Secondly, we define the function class $C_{f}^2(\mR^\infty)$ as follows: for any $\Phi\in C_{f}^2(\mR^\infty)$, there exist a $k\in\mN$ and a function $\phi\in C_b^2(\mR^k)$ such that $\Phi(\boldsymbol{x})=\phi(\boldsymbol{x}^1,\cdots,\boldsymbol{x}^k)$ for any 
$\boldsymbol{x}\in\mR^\infty$. That is, the functions in $C_{f}^2(\mR^\infty)$ only depend on the finite components of $\boldsymbol{x}$. Given $\Phi(\boldsymbol{x})=\phi(\boldsymbol{x}^1,\cdots,\boldsymbol{x}^k)\in C_{f}^2(\mR^\infty)$, we define $\partial_i\Phi(\boldsymbol{x}), \partial_{ij}\Phi(\boldsymbol{x}), i,j\geq 1$ by
\ce
\partial_i\Phi(\boldsymbol{x}):=\left\{\begin{array}{ll}
\partial_i\phi(\boldsymbol{x}^1,\cdots,\boldsymbol{x}^k), \quad i\leq k, \\
0, \qquad\qquad\qquad\quad i>k,
\end{array}
\right.
\quad 
\partial_{ij}\Phi(\boldsymbol{x}):=\left\{\begin{array}{ll}
\partial_{ij}\phi(\boldsymbol{x}^1,\cdots,\boldsymbol{x}^k), \quad i,j\leq k, \\
0, \qquad\qquad\qquad\qquad i,j>k.
\end{array}
\right.
\de

Fix two Borel measurable mappings
$$
\a: [0,T]\times\mR^\infty\mapsto\mS_+(\mR^{\infty}), \qquad \b: [0,T]\times\mR^\infty\mapsto\mR^\infty,
$$
where $\mS_+(\mR^{\infty})$ denotes the set of double sequences $(\a^{ij})$ satisfying 
$$
\sum\limits_{i,j=1}^k\a^{ij}\xi^i\xi^j\geq 0, \quad \forall k\geq 1, \quad \forall \xi\in\mR^\infty.
$$
And then we define the operator $\cL(\a,\b)$ on $C_{f}^2(\mR^\infty)$:
$$
\cL(\a,\b)\Phi(\boldsymbol{x}):=\frac{1}{2}\a^{ij}(t,\boldsymbol{x})\partial_{ij}\Phi(\boldsymbol{x})+\b^i(t,\boldsymbol{x})\partial_i\Phi(\boldsymbol{x}), \quad \Phi\in C_{f}^2(\mR^\infty).
$$
In the following we give out the definition of solutions to the martingale problem associated with $\cL(\a,\b)$.

\bd\label{martsolu}
For $\mu\in\cP(\mR^\infty)$. A probability measure $\boldsymbol{\eta}$ on $(C_T^\infty, \cB)$ is called a solution to the martingale problem associated with 
$\cL(\a,\b)$ with the initial law $\mu$ at time $0$, if

(i) $\boldsymbol{\eta}\circ e^{-1}_0=\mu$,

(ii)
\ce
\int_0^T\int_{\mR^\infty}\left(|\a^{ij}(s,\boldsymbol{x})|+|\b^i(s,\boldsymbol{x})|\right)\boldsymbol{\eta}\circ e_s^{-1}(\dif \boldsymbol{x})\dif s<\infty, \quad i,j\geq 1,
\de

(iii) For any $\Phi\in{C_f^2(\mR^\infty)}$,
\be
\cM_t^\Phi&:=&\Phi(w_t)-\Phi(w_0)-\int_0^t\cL(\a,\b)\Phi(w_s)ds
\label{eq2}
\ee
is a $\bar{\cB}_t$-adapted martingale under the probability measure $\boldsymbol{\eta}$. The uniqueness of solutions to the martingale problem associated with $\cL(\a,\b)$ with the initial law $\mu$ at time $0$ means that, if $\boldsymbol{\eta}, \tilde{\boldsymbol{\eta}}$ are two solutions to the martingale problem associated with $\cL(\a,\b)$ with the initial law $\mu$ at time $0$, then $\boldsymbol{\eta}\circ e_t^{-1}=\tilde{\boldsymbol{\eta}}\circ e_t^{-1}$ for any $t\in[0,T]$.
\ed

Consider the FPE associated with the operator $\cL(\a,\b)$:
\be
\partial_t\mQ_t=(\cL(\a,\b))^*\mQ_t,
\label{FPE1}
\ee
where $(\cL(\a,\b))^*$ is the adjoint operator of $\cL(\a,\b)$, and $(\mQ_t)_{t\in[0,T]}$ is a family of probability measures on $\sB(\mR^\infty)$.  Weak solutions of Eq.(\ref{FPE1}) are defined as follows. 

\bd\label{weakfpe}
A measurable family $(\mQ_t)_{t\in[0,T]}$ of probability measures on $\sB(\mR^\infty)$ is called a
 weak solution of the FPE (\ref{FPE1}) if 
\be
\int_0^T\int_{\mR^\infty}\left(|\a^{ij}(s,\boldsymbol{x})|+|\b^i(s,\boldsymbol{x})|\right)\mQ_s(\dif \boldsymbol{x})\dif s<\infty, \quad i,j\geq 1,
\label{deficond1}
\ee
and for all $\Phi\in C_f^2(\mR^\infty)$ and $t\in[0,T]$,
\be
<\mQ_t, \Phi>=<\mQ_0, \Phi>+\int_0^t<\mQ_s, \cL(\a,\b)\Phi>\dif s.
\label{deficond2}
\ee 
The uniqueness of the weak solutions to Eq.(\ref{FPE1}) means that, if $(\mQ_t)_{t\in[0,T]}$ and $(\tilde{\mQ}_t)_{t\in[0,T]}$ are two weak solutions to Eq.(\ref{FPE1}) with $\mQ_0=\tilde{\mQ}_0$, then $\mQ_t=\tilde{\mQ}_t$ for any $t\in[0,T]$.
\ed

Next, we describe the relationship between the solutions of the martingale problem associated with 
$\cL(\a,\b)$ and weak solutions of the FPE (\ref{FPE1}).

\bt(The superposition principle on $\mR^\infty$)\label{supeprinrn} For $\mu\in\cP(\mR^\infty)$.

(i) The existence of solutions $\boldsymbol{\eta}$ to the martingale problem associated with $\cL(\a,\b)$ with the initial law $\mu$ at time $0$ is equivalent to the existence of weak solutions $(\mQ_t)_{t\in[0,T]}$ of the FPE (\ref{FPE1}) with the initial value $\mu$. Moreover, $\boldsymbol{\eta}_t:=\boldsymbol{\eta}\circ e_t^{-1}=\mQ_t$ for any $t\in[0,T]$.

(ii) The uniqueness of solutions $\boldsymbol{\eta}$ to the martingale problem associated with $\cL(\a,\b)$ with the initial law $\mu$ at time $0$ is equivalent to the uniqueness of weak solutions $(\mQ_t)_{t\in[0,T]}$ of the FPE (\ref{FPE1}) with the initial value $\mu$.
\et
\begin{proof}
For (i), if $\boldsymbol{\eta}$ is a solution to the martingale problem associated with $\cL(\a,\b)$ with the initial law $\mu$ at time $0$, by (\ref{eq2}) in Definition \ref{martsolu} it is direct to obtain that $(\boldsymbol{\eta}_t)_{t\in[0,T]}$ is a weak solution of the FPE (\ref{FPE1}) with the initial value $\mu$. Conversely, if $(\mQ_t)_{t\in[0,T]}$ is a weak solution of the FPE (\ref{FPE1}) with the initial value $\mu$, we show that there exists a solution $\boldsymbol{\eta}$ to the martingale problem associated with $\cL(\a,\b)$ with the initial law $\mu$ at time $0$ satisfying $\boldsymbol{\eta}_t=\mQ_t$ for any $t\in[0,T]$. The proof is divided into three following steps.

{\bf Step 1.} We reduce the problem in $\mR^\infty$ to a sequence of problems in $\mR^k$ for any $k\in\mN$. 

First of all, set for any $k\in\mN$
$$
\pi^k: \mR^\infty\rightarrow \mR^k, \quad \pi^k(\boldsymbol{x})=(\boldsymbol{x}^1, \boldsymbol{x}^2, \cdots, \boldsymbol{x}^k), \quad \boldsymbol{x}\in\mR^\infty,
$$
and then $\tilde{\mQ}^k_t:=\mQ_t\circ(\pi^k)^{-1}\in\cP(\mR^k)$ for $t\in[0,T]$. Put
\ce
&&(\tilde{\a}^k)^{ij}(t,x):=\frac{\dif (\a^{ij}(t,\cdot)\mQ_t)\circ(\pi^k)^{-1}}{\dif \mQ_t\circ(\pi^k)^{-1}}(x),\\
&&(\tilde{\b}^k)^{i}(t,x):=\frac{\dif (\b^i(t,\cdot)\mQ_t)\circ(\pi^k)^{-1}}{\dif \mQ_t\circ(\pi^k)^{-1}}(x), \quad i,j=1, 2, \cdots, k, \quad x\in\mR^k.
\de
Note that
\be
&&\int_0^T\int_{\mR^k}\left(|(\tilde{\a}^k)^{ij}(s,x)|+|(\tilde{\b}^k)^{i}(s,x)|\right)\tilde{\mQ}^k_s(\dif x)\dif s\no\\
&\leq&\int_0^T\int_{\mR^\infty}\left(|\a^{ij}(s,\boldsymbol{x})|+|\b^i(s,\boldsymbol{x})|\right)\mQ_s(\dif \boldsymbol{x})\dif s<\infty,
\label{finiinfi}
\ee
and for any $\phi\in C_b^2(\mR^k)$ with $\Phi(\boldsymbol{x}):=\phi(\pi^k(\boldsymbol{x}))\in C_f^2(\mR^\infty)$, 
\ce
<\tilde{\mQ}^k_t, \phi>&=&<\mQ_t, \Phi>=<\mQ_0, \Phi>+\int_0^t<\mQ_s, \cL(\a,\b)\Phi>\dif s\\
&=&<\tilde{\mQ}^k_0, \phi>+\int_0^t<\tilde{\mQ}^k_s, \cL(\tilde{\a}^k, \tilde{\b}^k)\phi>\dif s,
\de
where $\cL(\tilde{\a}^k, \tilde{\b}^k)$ denotes the operator associated with $\tilde{\a}^k, \tilde{\b}^k$. So, by \cite[Definition 2.2]{T}, it holds that $(\tilde{\mQ}_t^k)_{t\in[0,T]}$ is a weak solution of the FPE associated with $\cL(\tilde{\a}^k, \tilde{\b}^k)$. And then by \cite[Theorem 2.5]{T}, we obtain that there exists a solution $\tilde{\boldsymbol{\eta}}^k\in\cP(C([0,T],\mR^k))$  to the martingale problem associated with 
$\cL(\tilde{\a}^k, \tilde{\b}^k)$ such that $\tilde{\boldsymbol{\eta}}^k_t=\tilde{\mQ}^k_t$.

{\bf Step 2.} Set 
\ce
&&\bar{\pi}^k: \mR^k\rightarrow \mR^\infty, \quad \bar{\pi}^k(x)=(x^1, x^2, \cdots, x^k, 0, 0, \cdots), \quad x\in\mR^k,\\
&&J^k: C([0,T], \mR^k)\rightarrow C([0,T], \mR^\infty), ~ J^k(\zeta)=(\zeta^1, \zeta^2, \cdots, \zeta^k, 0, 0, \cdots), ~\zeta\in C([0,T], \mR^k),
\de
and then $\mQ^k_t:=\tilde{\mQ}^k_t\circ(\bar{\pi}^k)^{-1}\in\cP(\mR^\infty)$ and $\boldsymbol{\eta}^k:=\tilde{\boldsymbol{\eta}}^k\circ(J^k)^{-1}\in\cP(C([0,T], \mR^\infty))$. Let 
\ce
(\a^k)^{ij}(t,\boldsymbol{x}):=\left\{\begin{array}{ll}
(\tilde{\a}^k)^{ij}(t,x), \quad i,j\leq k, \\
0, \qquad\qquad\quad i,j>k,
\end{array}
\right.
\quad
(\b^k)^i(t,\boldsymbol{x}):=\left\{\begin{array}{ll}
(\tilde{\b}^k)^i(t,x), \quad i\leq k, \\
0, \qquad\qquad\quad i>k,
\end{array}
\right.
\de
and $\cL(\a^k, \b^k)$ be the operator associated with $\a^k, \b^k$. Thus, the statement in {\bf Step 1.} yields that $\boldsymbol{\eta}^k$ is a solution to the martingale problem associated with $\cL(\a^k, \b^k)$ 
and satisfies $\boldsymbol{\eta}^k_t=\mQ^k_t$ for any $t\in[0,T]$. 

In the following, we prove that $\{\boldsymbol{\eta}^k\}$ is tight. By (\ref{deficond1}) and the de la Vall\'ee Poussin criterion, it holds that for $i\geq 1$
$$
\int_0^T\int_{\mR^\infty}\left[\Theta_1^i(|\b^i(t,\boldsymbol{x})|)+\Theta_2^i(|\a^{ii}(t,\boldsymbol{x})|)\right]\mQ_t(\dif\boldsymbol{x})\dif t<2^{-i}, 
$$
where $\Theta_1^i: [0,+\infty)\rightarrow[0,+\infty)$, $\Theta_2^i: [0,+\infty)\rightarrow[0,+\infty)$ are convex, lower semicontinuous and satisfy
$$
\lim\limits_{x\rightarrow+\infty}\frac{\Theta_1^i(x)}{x}=\lim\limits_{x\rightarrow+\infty}\frac{\Theta_2^i(x)}{x}=+\infty,
$$
and for some constant $C\geq 0$, $\Theta_2^i(2x)\leq \Theta_2^i(x)$ for any $x\geq 0$. Besides, we choose coercive functions $\t^i: [0,+\infty]\rightarrow[0,+\infty]$ with $\lim\limits_{x\rightarrow+\infty}\t^i(x)=+\infty$ such that 
$$
\int_{\mR^\infty}\t^i(|\boldsymbol{x}^i|)\mu(\dif\boldsymbol{x})<2^{-i}.
$$
Thus, from \cite[Corollary A5]{T}, it follows that there exist coercive functionals $\Psi^i: C([0,T], \mR)\rightarrow[0,\infty)$ associated with $\t^i, \Theta_1^i, \Theta_2^i$ such that 
\be
&&\int_{C_T^\infty}\Psi^i(\boldsymbol{x}^i\circ w)\dif\boldsymbol{\eta}^k(w)\no\\
&\leq& \int_{\mR^\infty}\t^i(|\boldsymbol{x}^i|)\dif\boldsymbol{\eta}_0^k(\boldsymbol{x})+\int_0^T\int_{\mR^\infty}\left[\Theta_1^i(|(\b^k)^i(t,\boldsymbol{x})|)
+\Theta_2^i(|(\a^k)^{ii}(t,\boldsymbol{x})|)\right]\dif\boldsymbol{\eta}_t^k(\boldsymbol{x})\dif t\no\\
&=& \int_{\mR^\infty}\t^i(|\boldsymbol{x}^i|)\dif\mQ_0^k(\boldsymbol{x})+\int_0^T\int_{\mR^\infty}\left[\Theta_1^i(|(\b^k)^i(t,\boldsymbol{x})|)
+\Theta_2^i(|(\a^k)^{ii}(t,\boldsymbol{x})|)\right]\dif\mQ_t^k(\boldsymbol{x})\dif t\no\\
&\leq&\int_{\mR^\infty}\t^i(|\boldsymbol{x}^i|)\dif\mu(\boldsymbol{x})+\int_0^T\int_{\mR^\infty}\left[\Theta_1^i(|\b^i(t,\boldsymbol{x})|)
+\Theta_2^i(|\a^{ii}(t,\boldsymbol{x})|)\right]\dif\mQ_t(\boldsymbol{x})\dif t\no\\
&\leq&2\cdot 2^{-i},
\label{coeresti}
\ee
where the second last inequality is baed on the Jensen inequality and (\ref{finiinfi}). For $w\in C_T^\infty$, set 
$$
\Psi(w):=\sum\limits_{i=1}^\infty\Psi^i(\boldsymbol{x}^i\circ w),
$$
and then the functional $\Psi$ is coercive and 
$$
\sup\limits_{k\in\mN}\int_{C_T^\infty}\Psi(w)\boldsymbol{\eta}^k(\dif w)\overset{(\ref{coeresti})}{<}\infty.
$$
From this, it follows that $\{\boldsymbol{\eta}^k\}$ is tight.

{\bf Step 3.} We prove that the limit point $\boldsymbol{\eta}$ of $\{\boldsymbol{\eta}^k\}$ is a solution of the martingale problem associated with 
$\cL(\a,\b)$ with the initial law $\mu$ at time $0$ such that $\boldsymbol{\eta}_t=\mQ_t$ for any $t\in[0,T]$.

First of all, note that $\boldsymbol{\eta}_t^k\circ(\pi^l)^{-1}=\mQ^k_t\circ(\pi^l)^{-1}$ for any $l\in\mN$. So, as $k\rightarrow\infty$, based on the continuity of two mpppings $e_t, \pi^l$, it holds that $\boldsymbol{\eta}_t\circ(\pi^l)^{-1}=\mQ_t\circ(\pi^l)^{-1}$. 
Since $\sB(\mR^\infty)$ is generated by the finite dimensional cylindrical sets, we know that $\boldsymbol{\eta}_t=\mQ_t$ for any $t\in[0,T]$.   

Next, we show that $\boldsymbol{\eta}$ is a solution of the martingale problem associated with 
$\cL(\a,\b)$ with the initial law $\mu$ at time $0$. That is, it is sufficient to check that for $0\leq s<t\leq T$ and a bounded continuous $\bar{\cB}_s$-measurable functional $\chi_s: C^\infty_T\mapsto \mR$,
\be
\int_{C^\infty_T}\left[\Phi(w_t)-\Phi(w_s)-\int_s^t(\cL(\a,\b)\Phi)(w_r)dr\right]\chi_s(w)\boldsymbol{\eta}(\dif w)=0, \quad \forall \Phi\in{C_f^2(\mR^\infty)}.
\label{esti01}
\ee

To prove (\ref{esti01}), we arbitrarily take $\Phi\in{C_f^2(\mR^\infty)}$ with $\Phi(\boldsymbol{x})=\phi(\pi^u(\boldsymbol{x}))$ for $u\in\mN$ and $\phi\in C^2_b(\mR^u)$. Besides, note that 
$$
\{\zeta\in L^1(\mQ_t(\dif \boldsymbol{x})\dif t):  \zeta(t, \cdot)\in C_f^2(\mR^\infty)~\mbox{for  all}~ t\in(0,T)\}
$$
is dense in $L^1(\mQ_t(\dif \boldsymbol{x})\dif t)$, and then
\ce
\bigcup_{i=1}^\infty&&\{\zeta\in L^1(\mQ_t(\dif \boldsymbol{x})\dif t):  \zeta(t, \cdot)\in C_f^2(\mR^\infty), \exists i\in \mN, \zeta(t, \boldsymbol{x})~\mbox{only depends on the front}~ \\ 
&&i ~\mbox{components of}~\boldsymbol{x},
 ~\mbox{for  all}~ t\in(0,T)\}
\de
is also dense in $L^1(\mQ_t(\dif \boldsymbol{x})\dif t)$. Thus by (\ref{deficond1}), we know that for any $\e>0$ and the coefficients $\b, \a$, there exist $\bar{\b}: [0,T]\times\mR^\infty\mapsto\mR^\infty, \bar{\a}: [0,T]\times\mR^\infty\mapsto \mS_+(\mR^\infty)$ such that 

$(i)$ $\bar{\b}(t,\cdot), \bar{\a}(t,\cdot)\in C_f^2(\mR^\infty)$ for all $t\in(0,T)$;

$(ii)$
$$
\int_0^T\int_{\mR^\infty}\left(\sum\limits_{i=1}^u|\b^i(t,\boldsymbol{x})-\bar{\b}^i(t,\boldsymbol{x})|+\sum\limits_{i,j=1}^u|\a^{ij}(t,\boldsymbol{x})-\bar{\a}^{ij}(t,\boldsymbol{x})|\right)\mQ_t(\dif \boldsymbol{x})\dif t<\e.
$$
And then the operators with respect to $\bar{\b}, \bar{\a}$ are denoted as $\cL(\bar{\a}, \bar{\b})$. 

Now, we treat (\ref{esti01}). Inserting $\cL(\bar{\a},\bar{\b})$, one can estimate (\ref{esti01}) to get
\be
&&\left|\int_{C^\infty_T}\left[\Phi(w_t)-\Phi(w_s)-\int_s^t(\cL(\a,\b)\Phi)(w_r)dr\right]\chi_s(w)\boldsymbol{\eta}(\dif w)\right|\no\\
&\leq&\left|\int_{C^\infty_T}\left[\Phi(w_t)-\Phi(w_s)-\int_s^t(\cL(\bar{\a},\bar{\b})\Phi)(w_r)dr\right]\chi_s(w)\boldsymbol{\eta}(\dif w)\right|\no\\
&&+\left|\int_{C^\infty_T}\left[\int_s^t\left((\cL(\bar{\a},\bar{\b})\Phi)(w_r)-(\cL(\a,\b)\Phi)(w_r)\right)dr\right]\chi_s(w)\boldsymbol{\eta}(\dif w)\right|\no\\
&=:&I_1+I_2.
\label{esti11}
\ee
To deal with $I_1$, we recall that $\boldsymbol{\eta}^k$ is a solution to the martingale problem associated with $\cL(\a^k, \b^k)$, which means that 
$$
\int_{C^\infty_T}\left[\Phi(w_t)-\Phi(w_s)-\int_s^t(\cL(\a^k, \b^k)\Phi)(w_r)dr\right]\chi_s(w)\boldsymbol{\eta}^k(\dif w)=0.
$$
So, by some simple calculation it holds that
\ce
&&\int_{C^\infty_T}\left[\Phi(w_t)-\Phi(w_s)-\int_s^t(\cL(\bar{\a},\bar{\b})\Phi)(w_r)dr\right]\chi_s(w)\boldsymbol{\eta}^k(\dif w)\\
&=&\int_{C^\infty_T}\left[\int_s^t\left((\cL(\a^k, \b^k)\Phi)(w_r)-(\cL(\bar{\a},\bar{\b})\Phi)(w_r)\right)dr\right]\chi_s(w)\boldsymbol{\eta}^k(\dif w),
\de
and furthermore
\be
&&\left|\int_{C^\infty_T}\left[\Phi(w_t)-\Phi(w_s)-\int_s^t(\cL(\bar{\a},\bar{\b})\Phi)(w_r)dr\right]\chi_s(w)\boldsymbol{\eta}^k(\dif w)\right|\no\\
&=&\left|\int_{C^\infty_T}\left[\int_s^t\left((\cL(\a^k, \b^k)\Phi)(w_r)-(\cL(\bar{\a},\bar{\b})\Phi)(w_r)\right)dr\right]\chi_s(w)\boldsymbol{\eta}^k(\dif w)\right|\no\\
&\leq&C\int_{C^\infty_T}\int_s^t\left|(\cL(\a^k, \b^k)\Phi)(w_r)-(\cL(\bar{\a},\bar{\b})\Phi)(w_r)\right|dr\boldsymbol{\eta}^k(\dif w)\no\\
&=&C\int_s^t\int_{\mR^\infty}\left|(\cL(\a^k, \b^k)\Phi)(\boldsymbol{x})-(\cL(\bar{\a},\bar{\b})\Phi)(\boldsymbol{x})\right|\mQ_r^k(\dif \boldsymbol{x})dr\no\\
&\leq&C\int_s^t\int_{\mR^\infty}\bigg[\left|((\b^k)^i(r,\boldsymbol{x})-\bar{\b}^i(r,\boldsymbol{x}))\partial_{i}\Phi(\boldsymbol{x})\right|\no\\
&&\qquad\qquad +\left|((\a^k)^{ij}(r,\boldsymbol{x})-\bar{\a}^{ij}(r,\boldsymbol{x}))\partial_{ij}\Phi(\boldsymbol{x})\right|\bigg]\mQ_r^k(\dif \boldsymbol{x})dr.
\label{00}
\ee
As $k\rightarrow\infty$, based on (i), $\boldsymbol{\eta}^k\overset{w.}{\rightarrow}\boldsymbol{\eta}$ and the fact that $\pi^k\rightarrow I$, $\bar{\pi}^k\rightarrow I$, where $I: \mR^\infty\rightarrow\mR^\infty$ is an identity mapping, we take the limit on two sides of the inequality (\ref{00}) and have that
\be
I_1&\leq&C\int_s^t\int_{\mR^\infty}\bigg[|(\b^i(r,\boldsymbol{x})-\bar{\b}^i(r,\boldsymbol{x}))\partial_{i}\Phi(\boldsymbol{x})|+\left|(\a^{ij}(r,\boldsymbol{x})-\bar{\a}^{ij}(r,\boldsymbol{x}))\partial_{ij}\Phi(\boldsymbol{x})\right|\bigg]\mQ_r(\dif \boldsymbol{x})dr\no\\
&\leq&C\int_0^T\int_{\mR^\infty}\left(\sum\limits_{i=1}^u|\b^i(r,\boldsymbol{x})-\bar{\b}^i(r,\boldsymbol{x})|+\sum\limits_{i,j=1}^u|\a^{ij}(r,\boldsymbol{x})-\bar{\a}^{ij}(r,\boldsymbol{x})|\right)\mQ_r(\dif \boldsymbol{x})\dif r\no\\
&\overset{(ii)}{<}&C\e.
\label{esti21}
\ee

In the following, we treat $I_2$. By the similar deduction to that for $I_1$, one can obtain that
\be
I_2&\leq&C\int_{C^\infty_T}\int_s^t\left|(\cL(\bar{\a},\bar{\b})\Phi)(w_r)-(\cL(\a,\b)\Phi)(w_r)\right|dr\boldsymbol{\eta}(\dif w)\no\\
&=&C\int_s^t\int_{\mR^\infty}\bigg[\Big|(\b^i(r,\boldsymbol{x})-\bar{\b}^i(r,\boldsymbol{x}))\partial_{i}\Phi(\boldsymbol{x})\Big|+\Big|(\a^{ij}(r,\boldsymbol{x})-\bar{\a}^{ij}(r,\boldsymbol{x}))\partial_{ij}\Phi(\boldsymbol{x})\Big|\bigg]\mQ_r(\dif \boldsymbol{x})dr\no\\
&\leq&C\int_0^T\int_{\mR^\infty}\left(\sum\limits_{i=1}^u|\b^i(r,\boldsymbol{x})-\bar{\b}^i(r,\boldsymbol{x})|+\sum\limits_{i,j=1}^u|\a^{ij}(r,\boldsymbol{x})-\bar{\a}^{ij}(r,\boldsymbol{x})|\right)\mQ_r(\dif \boldsymbol{x})\dif r\no\\
&\overset{(ii)}{<}&C\e.
\label{esti31}
\ee
Combining (\ref{esti21}) (\ref{esti31}) with (\ref{esti11}), we get that
\ce
\left|\int_{C^\infty_T}\left[\Phi(w_t)-\Phi(w_s)-\int_s^t(\cL(\a,\b)\Phi)(w_r)dr\right]\chi_s(w)\boldsymbol{\eta}(\dif w)\right|<C\e.
\de
Letting $\e\rightarrow 0$, we finally have (\ref{esti01}). 

Since the proof of (ii) heavily depends on (i) and is also  straight, we omit it. The proof is complete.
\end{proof}

\br
Theorem \ref{supeprinrn} generalizes \cite[Theorem 7.1]{ad}, \cite[Theorem 7.1]{T0} and \cite[Theorem 2.5]{T}.
\er

\section{The Zakai equations from nonlinear filtering problems}\label{Zakaiequ}

In this section, we introduce two types of Zakai equations from nonlinear filtering problems. And then we take all of them as stochastic differential equations (SDEs in short) and define their weak solutions. 

\subsection{The Zakai equations from nonlinear filtering problems with correlated noises}\label{Zakaiequ1}

In this section, we introduce the Zakai equations from nonlinear filtering problems with correlated noises and then define their weak solutions. 

\begin{enumerate}[\bf{Assumption:}]
\item
\end{enumerate}
\begin{enumerate}[($\mathbf{H}^1_{b_1, \sigma_0, \sigma_1}$)] 
\item For $t\in[0,T]$ and $x_1, x_2\in\mR^n$,
\ce
&|b_1(t,x_1)-b_1(t,x_2)|\leq L_1(t)|x_1-x_2|\kappa_1(|x_1-x_2|),\\
&\|\sigma_0(t,x_1)-\sigma_0(t,x_2)\|^2\leq L_1(t)|x_1-x_2|^{2}\kappa_2(|x_1-x_2|),\\
&\|\sigma_1(t,x_1)-\sigma_1(t,x_2)\|^2\leq L_1(t)|x_1-x_2|^{2}\kappa_3(|x_1-x_2|),
\de
where $|\cdot|$ and $\|\cdot\|$ denote the Hilbert-Schmidt norms of a vector and a matrix, respectively. Here $L_1(t)>0$ is an increasing function and $\kappa_i$ is a positive continuous
function, bounded on $[1,\infty)$ and satisfies
\ce
\lim\limits_{x\downarrow0}\frac{\kappa_i(x)}{\log x^{-1}}<\infty, \quad i=1, 2, 3.
\de
\end{enumerate}
\begin{enumerate}[($\mathbf{H}^2_{b_1, \sigma_0,\sigma_1}$)]
\item For $t\in[0,T]$ and $x\in\mR^n$,
$$
|b_1(t,x)|^2+\|\sigma_0(t,x)\|^2+\|\sigma_1(t,x)\|^2\leq K_1(t)(1+|x|)^2,
$$
where $K_1(t)>0$ is an increasing function.
\end{enumerate}
\begin{enumerate}[($\mathbf{H}^1_{b_2, \sigma_2}$)] 
\item For $t\in[0,T]$, $\sigma_2(t)$ is invertible, and
$$
|\sigma^{-1}_2(t)b_2(t,x)|\leq K_2,~{for}~{all}~t\in[0,T], x\in\mR^n,
$$
where $K_2>0$ is a constant.
\end{enumerate}

By Theorem 1.2 in \cite{q1}, the system (\ref{Eq1}) has a pathwise unique strong solution denoted as $(X_t,Y_t)$. Set
\ce
h(t,X_t):=\sigma_2^{-1}(t)b_2(t,X_t),
\de
\ce
\Gamma^{-1}_t:=\exp\bigg\{-\int_0^th^i(s,X_s)\dif W^i_s-\frac{1}{2}\int_0^t
\left|h(s,X_s)\right|^2\dif s\bigg\}.
\de
And then by ($\mathbf{H}^1_{b_2, \sigma_2}$), we know that 
\ce
\mE\left(\int_0^T
\left|h(s,X_s)\right|^2\dif s\right)<\infty,
\de
and then $\Gamma^{-1}_t$ is an exponential martingale. Define a measure $\tilde{\mP}$ via
$$
\frac{\dif \tilde{\mP}}{\dif \mP}=\Gamma^{-1}_T.
$$
By the Girsanov theorem for Brownian motions(e.g.Theorem 3.17 in \cite{jjas}), one can obtain that under the measure $\tilde{\mP}$, 
\be\label{tilw}
\tilde{W}_t:=W_t+\int_0^t h(s,X_s)\dif s
\ee
is an $(\mathscr{F}_t)$-adapted Brownian motion. Set
\ce
<\tilde{\mP}_t,\varphi>:=\tilde{\mE}[\varphi(X_t)\Gamma_t|\mathscr{F}_t^Y], \quad \varphi\in\cB_b(\mR^n),
\de
where $\tilde{\mE}$ denotes the expectation under the measure $\tilde{\mP}$ and $\mathscr{F}_t^Y$ is the $\sigma$-field generated by $\{Y_s, 0\leq s\leq t\}$. And then $\tilde{\mP}_t$ is called as the unnormalized filtering of $X_t$ with respect to 
$\mathscr{F}_t^Y$, and the equation satisfied by $\tilde{\mP}_t$ is called as the Zakai equation. Moreover, by \cite[Theorem 2.8]{q2}, we can obtain the following Zakai equation:
\be
<\tilde{\mP}_t,\varphi>&=&<\tilde{\mP}_0,\varphi>+\int_0^t<\tilde{\mP}_s,\cL_s \varphi>\dif s
+\int_0^t<\tilde{\mP}_s,\varphi h^l(s,\cdot)+\partial_i\varphi\sigma^{il}_1(s,\cdot)>\dif \tilde{W}^l_s,\no\\
&&\qquad\qquad\qquad\qquad \varphi\in\cC_c^\infty(\mR^n), \quad t\in[0,T],
\label{zakaieq0}
\ee
where the operater $\cL_s$ is defined as
\be
\cL_s:=b^i_1(s,\cdot)\partial_i+\frac{1}{2}\(\sigma_0\sigma_0^*(s,\cdot)\)^{ij}\partial_{ij}+\frac{1}{2}
\(\sigma_1\sigma_1^*(s,\cdot)\)^{ij}\partial_{ij}.
\label{sigope}
\ee

Next, we view the Zakai equation (\ref{zakaieq0}) as a stochastic differential equation and define its weak solution.

\bd\label{soluzakai}
$\{(\hat{\Omega}, \hat{\mathscr{F}}, \{\hat{\mathscr{F}}_t\}_{t\in[0,T]},\hat{\mP}), (\hat{\mu}_t,
\hat{W}_t)\}$ is called a weak solution of the Zakai equation (\ref{zakaieq0}), if the following holds:

(i) $(\hat{\Omega}, \hat{\mathscr{F}}, \{\hat{\mathscr{F}}_t\}_{t\in[0,T]},\hat{\mP})$ is a complete filtered
probability space;

(ii) $\hat{\mu}_t$ is a $\cM(\mR^n)$-valued $(\hat{\mathscr{F}}_t)$-adapted continuous process and $\hat{\mu}_0\in\cP(\mR^n)$;

(iii) $\hat{W}_t$ is an $m$-dimensional $(\hat{\mathscr{F}}_t)$-adapted Brownian motion;

(iv) For any $t\in[0,T]$,
$$
\hat{\mP}\left(\int_0^t\int_{\mR^n}\Big(|b_1(r,z)|+|h(r,z)|^2+\|\sigma_1(r,z)\|^2+\|\sigma_0\sigma_0^*(r,z)\|\Big)\hat{\mu}_r(\dif z)\dif r<\infty\right)=1;
$$

(v) $(\hat{\mu}_t, \hat{W}_t)$ satisfies the following equation
\be
<\hat{\mu}_t, \varphi>&=&<\hat{\mu}_0, \varphi>+\int_0^t<\hat{\mu}_s, \cL_s \varphi>\dif s
+\int_0^t<\hat{\mu}_s,\partial_i\varphi\sigma^{il}_1(s,\cdot)>\dif \hat{W}^l_s\no\\
&&+\int_0^t<\hat{\mu}_s, \varphi h^l(s,\cdot)>\dif \hat{W}^l_s, \quad \varphi\in C_c^\infty(\mR^n), \quad t\in[0,T].
\label{zakaieq2}
\ee
\ed

\br\label{aweaksoluzaka}
We remind that by the above deduction, $\{(\Omega, \mathscr{F}, \{\mathscr{F}_t\}_{t\in[0,T]}, \tilde{\mP}),
(\tilde{\mP}_t, \tilde{W}_t)\}$ is a weak solution of the Zakai equation (\ref{zakaieq0}) with the initial value $\tilde{\mP}_0$.
\er

\bd\label{launzakai}
The uniqueness of weak solutions for the Zakai equation (\ref{zakaieq0}) means that if there exist two weak solutions $\{(\hat{\Omega}^1, \hat{\mathscr{F}}^1, \{\hat{\mathscr{F}}^1_t\}_{t\in[0,T]}, \hat{\mP}^1), (\hat{\mu}^1_t,\hat{W}^1_t)\}$ and $\{(\hat{\Omega}^2, \hat{\mathscr{F}}^2, \{\hat{\mathscr{F}}^2_t\}_{t\in[0,T]}, \\ \hat{\mP}^2),(\hat{\mu}^2_t,\hat{W}^2_t)\}$ with $\hat{\mP}^1\circ(\hat{\mu}^1_0)^{-1}=\hat{\mP}^2\circ(\hat{\mu}^2_0)^{-1}$, then $\hat{\mP}^1\circ(\hat{\mu}^1_t)^{-1}=\hat{\mP}^2\circ(\hat{\mu}^2_t)^{-1}$ for any $t\in[0,T]$.
\ed

\subsection{The Zakai equations from nonlinear filtering problems with correlated sensor noises}\label{Zakaiequ2}

In this subsection, we introduce the Zakai equations from nonlinear filtering problems with correlated sensor noises and define their weak solutions.

We make the following hypotheses:
\begin{enumerate}[(i)]
\item $\check{b}_1, \check{\sigma}_1$ satisfy ($\mathbf{H}^1_{b_1, \sigma_0, \sigma_1}$)-($\mathbf{H}^2_{b_1,\sigma_0,\sigma_1}$), where $\check{b}_1, \check{\sigma}_1$ replace $b_1, \sigma_1$;
\item $\check{b}_2(t,x)$ is bounded for all $t\in[0,T], x\in\mR^n$;
\item $\check{\sigma}_2\check{\sigma}^*_2+\check{\sigma}_3\check{\sigma}^*_3=I_m,$ where $I_m$ is the $m$-order unit matrix.
\end{enumerate}

Under the above assumptions, by \cite[Theorem 1.2]{q1} there exists a unique strong solution of the system (\ref{Eq0}) denoted as $(\check{X}_t,\check{Y}_t)$. Put
\ce
&&V_t:=\check{\sigma}_2 W_t+\check{\sigma}_3 B_t, \\
&&\Xi^{-1}_t:=\exp\bigg\{-\int_0^t\check{b}_2^i(s,\check{X}_s)\dif V^i_s-\frac{1}{2}\int_0^t
\left|\check{b}_2(s,\check{X}_s)\right|^2\dif s\bigg\},
\de
and then we know that $\Xi^{-1}_t$ is an exponential martingale. In addition, define the probability measure 
$$
\frac{\dif \tilde{\check{\mP}}}{\dif \mP}:=\Xi^{-1}_T,
$$
and set
\ce
<\tilde{\check{\mP}}_t,\varphi>:=\tilde{\check{\mE}}[\varphi(\check{X}_t)\Xi_t|\mathscr{F}_t^{\check{Y}}], \qquad \varphi\in\cB_b(\mR^n),
\de
where $\tilde{\check{\mE}}$ stands for the expectation under the probability measure $\tilde{\check{\mP}}$. Thus, by \cite[Corollary 5.2]{q2}, the Zakai equation of the system (\ref{Eq0}) is given by
\be
<\tilde{\check{\mP}}_t, \varphi>&=&<\tilde{\check{\mP}}_0, \varphi>+\int_0^t<\tilde{\check{\mP}}_s, \check{\cL}_s \varphi>\dif s
+\int_0^t<\tilde{\check{\mP}}_s, \varphi\check{b}_2^l(s,\cdot)+\partial_i\varphi\(\check{\sigma}_1(s,\cdot)\check{\sigma}^*_2\)^{il}>\dif \tilde{V}^l_s,\no\\
&& \qquad\qquad \varphi\in C_c^\infty(\mR^n), \quad t\in[0,T],
\label{zakaieq01}
\ee
where $\tilde{V}_t:=V_t+\int_0^t \check{b}_2(s,\check{X}_s)\dif s$ and $\check{\cL_s}:=\check{b}^i_1(s,\cdot)\partial_i+\frac{1}{2}\(\check{\sigma}_1\check{\sigma}^*_1(s,\cdot)\)^{ij}\partial_{ij}$.

\br\label{120zakai}
If $\sigma_1=0$ in Eq.(\ref{zakaieq0}) and $\check{\sigma}_2=0$ in Eq.(\ref{zakaieq01}), both Eq.(\ref{zakaieq0}) and Eq.(\ref{zakaieq01}) become the usual Zakai equations (c.f. \cite[Theorem 3.3]{qd}, \cite{roge, blr}). That is, both Eq.(\ref{zakaieq0}) and Eq.(\ref{zakaieq01}) can include a lot of the Zakai equations.
\er

In the following, we take the Zakai equation (\ref{zakaieq01}) as a stochastic differential equation and define its weak solution.

\bd\label{soluzakai2}
$\{(\hat{\check{\Omega}}, \hat{\check{\sF}}, \{\hat{\check{\sF}}_t\}_{t\in[0,T]},\hat{\check{\mP}}), (\hat{\check{\mu}}_t,
\hat{\check{V}}_t)\}$ is called a weak solution of the Zakai equation (\ref{zakaieq01}), if the following holds:

(i) $(\hat{\check{\Omega}}, \hat{\check{\sF}}, \{\hat{\check{\sF}}_t\}_{t\in[0,T]},\hat{\check{\mP}})$ is a complete filtered
probability space;

(ii) $\hat{\check{\mu}}_t$ is a $\cM(\mR^n)$-valued $(\hat{\check{\sF}}_t)$-adapted continuous process and $\hat{\check{\mu}}_0\in\cP(\mR^n)$;

(iii) $(\hat{\check{V}}_t)$ is an $m$-dimensional $(\hat{\check{\sF}}_t)$-adapted Brownian motion;

(iv) For any $t\in[0,T]$,
$$
\hat{\check{\mP}}\left(\int_0^t\int_{\mR^n}\Big(|\check{b}_1(r,z)|+\|\check{\sigma}_1(r,z)\|^2+|\check{b}_2(r,z)|^2\Big)\hat{\check{\mu}}_r(\dif z)\dif r<\infty\right)=1;
$$

(v) $(\hat{\check{\mu}}_t, \hat{\check{V}}_t)$ satisfies the following equation
\be
<\hat{\check{\mu}}_t, \varphi>&=&<\hat{\check{\mu}}_t, \varphi>+\int_0^t<\hat{\check{\mu}}_s, \check{\cL}_s \varphi>\dif s
+\int_0^t<\hat{\check{\mu}}_s, \varphi\check{b}_2^l(s,\cdot)+\partial_i\varphi\check{\sigma}^{ik}_1(s,\cdot)\check{\sigma}^{lk}_2>\dif \hat{\check{V}}^l_s, \no\\
&& \qquad\qquad \quad \varphi\in C_c^\infty(\mR^n), \quad t\in[0,T].
\label{zakaieq22}
\ee
\ed

\br\label{aweaksoluzaka2}
By the above deduction, we have that $\{(\Omega, \mathscr{F}, \{\mathscr{F}_t\}_{t\in[0,T]}, \tilde{\check{\mP}}),
(\tilde{\check{\mP}}_t, \tilde{V}_t)\}$ is a weak solution of the Zakai equation (\ref{zakaieq01}) with the initial value $\tilde{\check{\mP}}_0$.
\er

\bd\label{launzakai}
The uniqueness of weak solutions for the Zakai equation (\ref{zakaieq01}) means that if there exist two weak solutions $\{(\hat{\check{\Omega}}^1, \hat{\check{\mathscr{F}}}^1, \{\hat{\check{\mathscr{F}}}^1_t\}_{t\in[0,T]}, \hat{\check{\mP}}^1), (\hat{\check{\mu}}^1_t,\hat{\check{W}}^1_t)\}$ and $\{(\hat{\check{\Omega}}^2, \hat{\check{\mathscr{F}}}^2, \{\hat{\check{\mathscr{F}}}^2_t\}_{t\in[0,T]}, \\\hat{\check{\mP}}^2),(\hat{\check{\mu}}^2_t,
\hat{\check{W}}^2_t)\}$ with $\hat{\check{\mP}}^1\circ(\hat{\check{\mu}}^1_0)^{-1}=\hat{\check{\mP}}^2\circ(\hat{\check{\mu}}^2_0)^{-1}$, then
$\hat{\check{\mP}}^1\circ(\hat{\check{\mu}}^1_t)^{-1}=\hat{\check{\mP}}^2\circ(\hat{\check{\mu}}^2_t)^{-1}$ for any $t\in[0,T]$.
\ed

\section{The Fokker-Planck equations on $\cM(\mR^n)$}\label{TFPE}

In the section, we introduce two types of Fokker-Planck equations associated with the Zakai equations (\ref{zakaieq0}) and (\ref{zakaieq01}).

\subsection{A type of Fokker-Planck equations associated with the Zakai equations (\ref{zakaieq0})}\label{fpes1}

In the subsection, we introduce a type of Fokker-Planck equations associated with the Zakai equations (\ref{zakaieq0}) and define their weak solutions.

To do this, we introduce the following function class $\sS$ on $\cM(\mR^n)$:
$$
\sS:=\Big\{\mu\mapsto g\left(<\mu,\varphi_1>,\cdots,<\mu,\varphi_k>\right): k\in\mN, g\in C_b^2(\mR^{k}), \varphi_1, \cdots, \varphi_k\in C_c^\infty(\mR^n)\Big\}.
$$
And then for any $G(\mu)=g\left(<\mu,\varphi_1>,\cdots,<\mu,\varphi_k>\right)=:g(<\mu,\boldsymbol{\varphi}>)\in\sS$, we define an operator ${\bf L}_t$ on $\sS$ as follows:
\ce
{\bf L}_tG(\mu)&=&\frac{1}{2}\partial_{uv} g(<\mu,\boldsymbol{\varphi}>)<\mu,\varphi_u h^l(t,\cdot)+\partial_i\varphi_u\sigma^{il}_1(t,\cdot)><\mu,\varphi_v h^l(t,\cdot)+\partial_i\varphi_v\sigma^{il}_1(t,\cdot)>\\
&&+\partial_u g(<\mu,\boldsymbol{ \varphi}>)<\mu,\cL_t\varphi_u>\\
&\overset{(\ref{sigope})}{=}&\frac{1}{2}\partial_{uv} g(<\mu,\boldsymbol{\varphi}>)<\mu,\varphi_u h^l(t,\cdot)+\partial_i\varphi_u\sigma^{il}_1(t,\cdot)><\mu,\varphi_v h^l(t,\cdot)+\partial_i\varphi_v\sigma^{il}_1(t,\cdot)>\\
&&+\partial_u g(<\mu,\boldsymbol{ \varphi}>)<\mu,\partial_i\varphi_u b^i_1(t,\cdot)>+\frac{1}{2}\partial_u g(<\mu,\boldsymbol{ \varphi}>)<\mu,\partial_{ij}\varphi_u\left(\sigma_0\sigma_0^*(t,\cdot)\right)^{ij}>\\
&&+\frac{1}{2}\partial_u g(<\mu,\boldsymbol{ \varphi}>)<\mu,\partial_{ij}\varphi_u\left(\sigma_1\sigma_1^*(t,\cdot)\right)^{ij}>.
\de

\br
 To understand the operator ${\bf L}_t$, we explain it by a special case. Set $b_2(t,x)=0$, and then $h(t,x)=0$ and ${\bf L}_t$ becomes
 \be
{\bf L}_tG(\mu)&=&\frac{1}{2}\partial_{uv} g(<\mu,\boldsymbol{\varphi}>)<\mu,\partial_i\varphi_u\sigma^{il}_1(t,\cdot)><\mu,\partial_i\varphi_v\sigma^{il}_1(t,\cdot)>\no\\
&&+\frac{1}{2}\partial_u g(<\mu,\boldsymbol{ \varphi}>)<\mu,\partial_{ij}\varphi_u\left(\sigma_0(t,\cdot)\sigma_0^*(t,\cdot)\right)^{ij}>\no\\
&&+\frac{1}{2}\partial_u g(<\mu,\boldsymbol{ \varphi}>)<\mu,\partial_{ij}\varphi_u\left(\sigma_1(t,\cdot)\sigma_1^*(t,\cdot)\right)^{ij}>\no\\
&&+\partial_u g(<\mu,\boldsymbol{ \varphi}>)<\mu,\partial_i\varphi_u b^i_1(t,\cdot)>.
\label{hequ0}
\ee
For $\mu\in\cP_2(\mR^n)$, based on the definition of L-derivatives in Subsection \ref{lde}, (\ref{hequ0}) goes into
 \be
{\bf L}_tG(\mu)&\overset{(\ref{cylilderi})}{=}&\frac{1}{2}\int_{\mR^n}\int_{\mR^n}\tr\(\partial^2_{\mu}G(\mu)(y,y')\sigma_1(t,y)\sigma^*_1(t,y')\)\mu(\dif y)\mu(\dif y')\no\\
&&+\frac{1}{2}\int_{\mR^n}\tr\(\partial_y\partial_{\mu}G(\mu)(y)\sigma_0(t,y)\sigma^*_0(t,y)\)\mu(\dif y)\no\\
&&+\frac{1}{2}\int_{\mR^n}\tr\(\partial_y\partial_{\mu}G(\mu)(y)\sigma_1(t,y)\sigma^*_1(t,y)\)\mu(\dif y)\no\\
&&+\int_{\mR^n}\partial_{\mu}G(\mu)(y)b_1(t,y)\mu(\dif y).
\label{hequ01}
\ee
Thus, the operator ${\bf L}_t$ is just right the integral of the L-differential operators.
\er

Consider the following FPE:
\be
\partial_t\Lambda_t={\bf L}^*_t\Lambda_t,
\label{fpecoup}
\ee
where $(\Lambda_t)_{t\in[0,T]}$ is a family of probability measures on $\sB(\cM(\mR^n))$. And then weak solutions of the FPE (\ref{fpecoup}) are defined as follows.

\bd\label{weaksolufpe}
A measurable family $(\Lambda_t)_{t\in[0,T]}$ of probability measures on $\sB(\cM(\mR^n))$ is called a weak solution of the FPE (\ref{fpecoup}) if 
\be
\int_0^T\int_{\cM(\mR^n)}\int_{\mR^n}\Big(|b_1(r,z)|+|h(r,z)|^2+\|\sigma_1(r,z)\|^2+\|\sigma_0\sigma_0^*(r,z)\|\Big)\mu(\dif z)\Lambda_r(\dif \mu)\dif r<\infty, \label{fpem01}
\ee
and for any $G\in\sS$ and $0\leq t\leq T$,
\be
\int_{\cM(\mR^n)}G(\mu)\Lambda_t(\dif\mu)&=&\int_{\cM(\mR^n)}G(\mu)\Lambda_0(\dif \mu)+\int_0^t\int_{\cM(\mR^n)}{\bf L}_rG(\mu)\Lambda_r(\dif \mu)\dif r.
\label{fpem2}
\ee
The uniqueness of the weak solutions to Eq.(\ref{fpecoup}) means that, if $(\Lambda_t)_{t\in[0,T]}$ and $(\tilde{\Lambda}_t)_{t\in[0,T]}$ are two weak solutions to Eq.(\ref{fpecoup}) with $\Lambda_0=\tilde{\Lambda}_0$, then $\Lambda_t=\tilde{\Lambda}_t$ for any $t\in[0,T]$.
\ed

It is easy to see that under the condition (\ref{fpem01}), the integral in the right side of Eq.(\ref{fpem2}) is well defined.

\subsection{A type of Fokker-Planck equations associated with the Zakai equations (\ref{zakaieq01})}\label{fpes2}

In the subsection, we introduce the other type of Fokker-Planck equations associated with the Zakai equations (\ref{zakaieq01}) and their weak solutions.

First of all, we define an operator ${\bf\check{L}}_t$ on $\sS$ as follows: for $G(\mu)=g(<\mu,\boldsymbol{\varphi}>)\in\sS$
\ce
{\bf\check{L}}_tG(\mu)&=&\frac{1}{2}\partial_{uv} g(<\mu,\boldsymbol{ \varphi}>)<\mu,\varphi_u\check{b}_2^l(t,\cdot)+\partial_i\varphi_u\check{\sigma}^{ik}_1(t,\cdot)\check{\sigma}^{lk}_2>\\
&&\times<\mu,\varphi_v\check{b}_2^l(t,\cdot)+\partial_i\varphi_v\check{\sigma}^{ik}_1(t,\cdot)\check{\sigma}^{lk}_2>+\partial_u g(<\mu,\boldsymbol{ \varphi}>)<\mu,\check{\cL}_t\varphi_u>\\
&=&\frac{1}{2}\partial_{uv} g(<\mu,\boldsymbol{ \varphi}>)<\mu,\varphi_u\check{b}_2^l(t,\cdot)+\partial_i\varphi_u\check{\sigma}^{ik}_1(t,\cdot)\check{\sigma}^{lk}_2>\\
&&\times<\mu,\varphi_v\check{b}_2^l(t,\cdot)+\partial_i\varphi_v\check{\sigma}^{ik}_1(t,\cdot)\check{\sigma}^{lk}_2>+\partial_u g(<\mu,\boldsymbol{ \varphi}>)<\mu,\partial_i\varphi_u \check{b}^i_1(t,\cdot)>\\
&&+\frac{1}{2}\partial_u g(<\mu,\boldsymbol{ \varphi}>)<\mu,\partial_{ij}\varphi_u\left(\check{\sigma}_1\check{\sigma}_1^*(t,\cdot)\right)^{ij}>.
\de
Specially, if $\check{b}_2(t,x)=0$, by the definition of L-derivatives in Subsection \ref{lde} it holds that for $\mu\in\cP_2(\mR^n)$
 \ce
 {\bf\check{L}}_tG(\mu)&=&\frac{1}{2}\partial_{uv} g(<\mu,\boldsymbol{ \varphi}>)<\mu,\partial_i\varphi_u\check{\sigma}^{ik}_1(t,\cdot)\check{\sigma}^{lk}_2><\mu,\partial_i\varphi_v\check{\sigma}^{ik}_1(t,\cdot)\check{\sigma}^{lk}_2>\\
 &&+\partial_u g(<\mu,\boldsymbol{ \varphi}>)<\mu,\partial_i\varphi_u \check{b}^i_1(t,\cdot)>+\frac{1}{2}\partial_u g(<\mu,\boldsymbol{ \varphi}>)<\mu,\partial_{ij}\varphi_u\left(\check{\sigma}_1\check{\sigma}_1^*(t,\cdot)\right)^{ij}>\\
 &=&\frac{1}{2}\int_{\mR^n}\int_{\mR^n}\tr\(\partial^2_{\mu}G(\mu)(y,y')\check{\sigma}_1(t,y)\check{\sigma}_2\check{\sigma}^*_2\sigma^*_1(t,y')\)\mu(\dif y)\mu(\dif y')\no\\
&&+\int_{\mR^n}\partial_{\mu}G(\mu)(y)\check{b}_1(t,y)\mu(\dif y)+\frac{1}{2}\int_{\mR^n}\tr\(\partial_y\partial_{\mu}G(\mu)(y)\check{\sigma}_1(t,y)\check{\sigma}^*_1(t,y)\)\mu(\dif y).
\de
That is, the operator ${\bf\check{L}}_t$ is also the integral of the L-differential operators.

Next, consider the following FPE:
\be
\partial_t\check{\Lambda}_t={\bf\check{L}}^*_t\check{\Lambda}_t,
\label{fpecoup2}
\ee
where $(\check{\Lambda}_t)_{t\in[0,T]}$ is a family of probability measures on $\sB(\cM(\mR^n))$. And then weak solutions of the FPE (\ref{fpecoup2}) are defined as follows.

\bd\label{weaksolufpe2}
A measurable family $(\check{\Lambda}_t)_{t\in[0,T]}$ of probability measures on $\sB(\cM(\mR^n))$ is called a weak solution of the FPE (\ref{fpecoup2}) if 
\be
\int_0^T\int_{\cM(\mR^n)}\int_{\mR^n}\Big(|\check{b}_1(r,z)|+\|\check{\sigma}_1(r,z)\|^2+|\check{b}_2(r,z)|^2\Big)\mu(\dif z)\check{\Lambda}_r(\dif \mu)\dif r<\infty, \label{fpem012}
\ee
and for any $G\in\sS$ and $0\leq t\leq T$,
\be
\int_{\cM(\mR^n)}G(\mu)\check{\Lambda}_t(\dif\mu)&=&\int_{\cM(\mR^n)}G(\mu)\check{\Lambda}_0(\dif \mu)+\int_0^t\int_{\cM(\mR^n)}\check{{\bf L}}_rG(\mu)\check{\Lambda}_r(\dif \mu)\dif r.
\label{fpem22}
\ee
The uniqueness of the weak solutions to Eq.(\ref{fpecoup}) means that, if $(\check{\Lambda}_t)_{t\in[0,T]}$ and $(\tilde{\check{\Lambda}}_t)_{t\in[0,T]}$ are two weak solutions to Eq.(\ref{fpecoup}) with $\check{\Lambda}_0=\tilde{\check{\Lambda}}_0$, then $\check{\Lambda}_t=\tilde{\check{\Lambda}}_t$ for any $t\in[0,T]$.
\ed

It is easy to see that under the condition (\ref{fpem012}), the integral in the right side of Eq.(\ref{fpem22}) is well defined.

\br\label{120fpe}
If $\sigma_1=0$ in Eq.(\ref{fpecoup}) and $\check{\sigma}_2=0$ in Eq.(\ref{fpecoup2}), both Eq.(\ref{fpecoup}) and Eq.(\ref{fpecoup2}) are just those FPEs associated with the usual Zakai equations. 
\er

\section{The superposition principles for the Zakai equations and Fokker-Plank equations on $\cM(\mR^n)$}\label{susupeprin}

In the section, we state and prove two superposition principles for the Zakai equations and Fokker-Plank equations on $\cM(\mR^n)$.

\subsection{The superposition principle between Eq.(\ref{zakaieq0}) and Eq.(\ref{fpecoup})}\label{susupeprin1}

In the subsection, we show the superposition between Eq.(\ref{zakaieq0}) and Eq.(\ref{fpecoup}). Let us begin with Eq.(\ref{zakaieq0}).

\bp\label{sdefpe}
Assume that $\{(\hat{\Omega}, \hat{\mathscr{F}}, \{\hat{\mathscr{F}}_t\}_{t\in[0,T]},\hat{\mP}), (\hat{\mu}_t,\hat{W}_t)\}$ is a weak solution of the Zakai equation (\ref{zakaieq0}). Then $(\sL_{\hat{\mu}_t})_{t\in[0,T]}$ is a weak solution of Eq.(\ref{fpecoup}), where $\sL_{\hat{\mu}_t}$ stands for the distribution of $\hat{\mu}_t$.
\ep
\begin{proof}
First of all, we justify that the condition (\ref{fpem01}) holds for $\Lambda_t=\sL_{\hat{\mu}_t}$. Since $\{(\hat{\Omega}, \hat{\mathscr{F}}, \\\{\hat{\mathscr{F}}_t\}_{t\in[0,T]},\hat{\mP}), (\hat{\mu}_t,\hat{W}_t)\}$ is a weak solution of the Zakai equation (\ref{zakaieq0}), it holds by (iv) in Definition \ref{soluzakai}  that
\ce
&&\int_0^T\int_{\cM(\mR^n)}\int_{\mR^n}\Big(|b_1(r,z)|+|h(r,z)|^2+\|\sigma_1(r,z)\|^2+\|\sigma_0\sigma_0^*(r,z)\|\Big)\mu(\dif z)\Lambda_r(\dif \mu)\dif r\\
&=&\int_0^T\hat{\mE}\int_{\mR^n}\Big(|b_1(r,z)|+|h(r,z)|^2+\|\sigma_1(r,z)\|^2+\|\sigma_0\sigma_0^*(r,z)\|\Big)\hat{\mu}_r(\dif z)\dif r\\
&=&\hat{\mE}\int_0^T\int_{\mR^n}\Big(|b_1(r,z)|+|h(r,z)|^2+\|\sigma_1(r,z)\|^2+\|\sigma_0\sigma_0^*(r,z)\|\Big)\hat{\mu}_r(\dif z)\dif r\\
&<&\infty,
\de
where $\hat{\mE}$ denotes the expectation under the probability measure $\hat{\mP}$.

Next, we deal with the verification of Eq.(\ref{fpem2}). By (\ref{zakaieq2}) in Definition \ref{soluzakai}, we know that $\hat{\mu}_t$ satisfies the following equation
\ce
<\hat{\mu}_t, \varphi_u>&=&<\hat{\mu}_0, \varphi_u>+\int_0^t<\hat{\mu}_s, \cL_s \varphi_u>\dif s
+\int_0^t<\hat{\mu}_s,\partial_i\varphi_u\sigma^{il}_1(s,\cdot)>\dif \hat{W}^l_s\no\\
&&+\int_0^t<\hat{\mu}_s, \varphi_u h^l(s,\cdot)>\dif \hat{W}^l_s, \quad \varphi_u\in \cC_c^\infty(\mR^n), u=1,2,\cdots,k.
\de
Thus, for any $G(\mu)=g(<\mu,\boldsymbol{\varphi}>)\in\sS$, applying the It\^o formula to the process $G(\hat{\mu}_t)$ and taking the expectation on two sides under the probability measure $\hat{\mP}$, one can obtain that
\ce
\hat{\mE}G(\hat{\mu}_t)&=&\hat{\mE}G(\hat{\mu}_0)+\int_0^t\hat{\mE}\partial_u g(<\hat{\mu}_s,\boldsymbol{ \varphi}>)<\hat{\mu}_s,\cL_s\varphi_u>\dif s\\
&&+\frac{1}{2}\int_0^t\hat{\mE}\partial_{uv} g(<\hat{\mu}_s,\boldsymbol{\varphi}>)<\hat{\mu}_s,\varphi_u h^l(s,\cdot)+\partial_i\varphi_u\sigma^{il}_1(s,\cdot)>\\
&&\qquad \times<\hat{\mu}_s,\varphi_v h^l(s,\cdot)+\partial_i\varphi_v\sigma^{il}_1(s,\cdot)>\dif s,
\de
and furthermore
\ce
\int_{\cM(\mR^n)}G(\mu)\Lambda_t(\dif\mu)&=&\int_{\cM(\mR^n)}G(\mu)\Lambda_0(\dif \mu)+\int_0^t\int_{\cM(\mR^n)}{\bf L}_rG(\mu)\Lambda_r(\dif \mu)\dif r.
\de
Thus, $\Lambda_t$ satisfies Eq.(\ref{fpem2}). The proof is complete.
\end{proof}

By Remark \ref{aweaksoluzaka} and Proposition \ref{sdefpe}, we obtain the following conclusion.

\bc\label{fpeweakexis}
Assuem that ($\mathbf{H}^1_{b_1, \sigma_0, \sigma_1}$) ($\mathbf{H}^2_{b_1, \sigma_0, \sigma_1}$) ($\mathbf{H}^1_{b_2, \sigma_2}$) hold. Then Eq.(\ref{fpecoup}) has a weak solution.
\ec

\bp\label{supeprinpn}
Suppose that $(\Lambda_t)_{t\in[0,T]}$ is a weak solution of Eq.(\ref{fpecoup}). Then there exists a weak solution $\{(\hat{\Omega}, \hat{\mathscr{F}}, \{\hat{\mathscr{F}}_t\}_{t\in[0,T]},\hat{\mP}), (\hat{\mu}_t,\hat{W}_t)\}$ of Eq.(\ref{zakaieq0}) 
such that $\Lambda_t=\sL_{\hat{\mu}_t}$.
\ep
\begin{proof}
{\bf Step 1.} We project Eq.(\ref{fpecoup}) to a corresponding equation on $\mR^\infty$ and prove that Eq.(\ref{zakaieq0}) has a weak solution on $\mR^\infty$.

First of all, for any $\{\varphi_u, u\in\mN\}\subset C_c^\infty(\mR^n)$, set
$$
\cT: \cM(\mR^n)\rightarrow\mR^\infty, \quad \cT(\mu)=(<\mu,\varphi_1>,\cdots,<\mu,\varphi_u>,\cdots),
$$
and then it holds that for any $\Phi(\boldsymbol{x})=\phi(\pi^k(\boldsymbol{x}))\in C_{f}^2(\mR^\infty)$ and $\mu\in\cM(\mR^n)$,
$$
\Phi(\cT(\mu))=\Phi(<\mu,\varphi_1>,\cdots,<\mu,\varphi_u>,\cdots)=\phi(<\mu,\varphi_1>,\cdots,<\mu,\varphi_k>),
$$
and $\Phi\circ\cT\in\sS$. Since $(\Lambda_t)_{t\in[0,T]}$ is a weak solution of Eq.(\ref{fpecoup}), it follows from Eq.(\ref{fpem2}) that for $\Phi\circ\cT\in\sS$
\ce
\int_{\cM(\mR^n)}(\Phi\circ\cT)(\mu)\Lambda_t(\dif\mu)&=&\int_{\cM(\mR^n)}(\Phi\circ\cT)(\mu)\Lambda_0(\dif \mu)+\int_0^t\int_{\cM(\mR^n)}{\bf L}_r(\Phi\circ\cT)(\mu)\Lambda_r(\dif \mu)\dif r,
\de
that is,
\ce
\int_{\mR^\infty}\Phi(\boldsymbol{x})\mQ_t(\dif\boldsymbol{x})&=&\int_{\mR^\infty}\Phi(\boldsymbol{x})\mQ_0(\dif\boldsymbol{x})+\int_0^t\int_{\mR^\infty}\sum\limits_{u=1}^\infty\partial_u\Phi(\boldsymbol{x})<\cT^{-1}(\boldsymbol{x}),\cL_r\varphi_u>\mQ_r(\dif\boldsymbol{x})\dif r\\
&&+\frac{1}{2}\int_0^t\int_{\mR^\infty}\sum\limits_{u=1}^\infty\sum\limits_{v=1}^\infty\partial_{uv}\Phi(\boldsymbol{x})<\cT^{-1}(\boldsymbol{x}),\varphi_u h^l(r,\cdot)+\partial_i\varphi_u\sigma^{il}_1(r,\cdot)>\\
&&\qquad\qquad \times<\cT^{-1}(\boldsymbol{x}),\varphi_v h^l(r,\cdot)+\partial_i\varphi_v\sigma^{il}_1(r,\cdot)>\mQ_r(\dif\boldsymbol{x})\dif r,
\de
where $\mQ_t:=\Lambda_t\circ\cT^{-1}$. Here we remind that since $\Phi(\boldsymbol{x})$ only depends on the finite components of $\boldsymbol{x}$, the infinite summations in the right side of the above equation are in fact finite ones.

Next, set for $u,v\in\mN$
\ce
&&\beta^u(r,\boldsymbol{x}):=<\cT^{-1}(\boldsymbol{x}),\cL_r\varphi_u>,\\
&&\a^{uv}(r,\boldsymbol{x}):=<\cT^{-1}(\boldsymbol{x}),\varphi_u h^l(r,\cdot)+\partial_i\varphi_u\sigma^{il}_1(r,\cdot)>\\
&&\qquad\qquad\qquad\times <\cT^{-1}(\boldsymbol{x}),\varphi_v h^l(r,\cdot)+\partial_i\varphi_v\sigma^{il}_1(r,\cdot)>,
\de
and then it holds that
\ce
&&\int_0^T\int_{\mR^\infty}|\b^u(r,\boldsymbol{x})|\mQ_r(\dif \boldsymbol{x})\dif r\\
&=&\int_0^T\int_{\mR^\infty}|<\cT^{-1}(\boldsymbol{x}),\cL_r\varphi_u>|\mQ_r(\dif \boldsymbol{x})\dif r=\int_0^T\int_{\cM(\mR^n)}|<\mu,\cL_r\varphi_u>|\Lambda_r(\dif\mu)\dif r\\
&\leq&\int_0^T\int_{\cM(\mR^n)}\int_{\mR^n}\bigg[\left|\partial_i\varphi(x)b_1^i(r,x)\right|+\frac{1}{2}\left|\partial_{ij}\varphi(x)(\sigma_0\sigma^*_0)^{ij}(r,x)\right|\\
&&\qquad\qquad +\frac{1}{2}\left|\partial_{ij}\varphi(x)(\sigma_1\sigma^*_1)^{ij}(r,x)\right|\bigg]\mu(\dif x)\Lambda_r(\dif\mu)\dif r\\
&\leq&C\int_0^T\int_{\cM(\mR^n)}\int_{\mR^n}\Big(|b_1(r,x)|+\|\sigma_0\sigma_0^*(r,x)\|+\|\sigma_1(r,x)\|^2\Big)\mu(\dif x)\Lambda_r(\dif \mu)\dif r\\
&\overset{(\ref{fpem01})}{<}&\infty,
\de
and
\ce
&&\int_0^T\int_{\mR^\infty}|\a^{uv}(r,\boldsymbol{x})|\mQ_r(\dif \boldsymbol{x})\dif r\\
&=&\int_0^T\int_{\mR^\infty}|<\cT^{-1}(\boldsymbol{x}),\varphi_u h^l(r,\cdot)+\partial_i\varphi_u\sigma^{il}_1(r,\cdot)>|\\
&&\qquad\qquad \times|<\cT^{-1}(\boldsymbol{x}),\varphi_v h^l(r,\cdot)+\partial_i\varphi_v\sigma^{il}_1(r,\cdot)>|\mQ_r(\dif \boldsymbol{x})\dif r\\
&=&\int_0^T\int_{\cM(\mR^n)}|<\mu,\varphi_u h^l(r,\cdot)+\partial_i\varphi_u\sigma^{il}_1(r,\cdot)>|\\
&&\qquad\qquad \times|<\mu,\varphi_v h^l(r,\cdot)+\partial_i\varphi_v\sigma^{il}_1(r,\cdot)>|\Lambda_r(\dif \mu)\dif r\\
&\leq&\int_0^T\int_{\cM(\mR^n)}\int_{\mR^n}(|\varphi_u(x)h^l(r,x)|+|\partial_i\varphi_u(x)\sigma^{il}_1(r,x)|)\mu(\dif x)\\
&&\qquad\qquad \times\int_{\mR^n}(|\varphi_v(x)h^l(r,x)|+|\partial_i\varphi_v(x)\sigma^{il}_1(r,x)|)\mu(\dif x)\Lambda_r(\dif \mu)\dif r\\
&\leq&C\int_0^T\int_{\cM(\mR^n)}\left(\int_{\mR^n}(|h(r,x)|+\|\sigma_1(r,x)\|)\mu(\dif x)\right)^2\Lambda_r(\dif \mu)\dif r\\
&\leq&C\int_0^T\int_{\cM(\mR^n)}\int_{\mR^n}(|h(r,x)|^2+\|\sigma_1(r,x)\|^2)\mu(\dif x)\Lambda_r(\dif \mu)\dif r\\
&\overset{(\ref{fpem01})}{<}&\infty.
\de
Thus, by Theorem \ref{supeprinrn}, we know that there exists a solution $\boldsymbol{\eta}$  to the martingale problem associated with 
$\cL(\a,\b)$ with the initial law $\mQ_0$ at time $0$ such that $\boldsymbol{\eta}_t=\mQ_t$ for any $t\in[0,T]$. By the similar deduction to that in \cite[Proposition 4.6]{ks}, it holds that there is an 
$m$-dimensional Brownian motion $\hat{W}$ defined on an extension $(\hat{\Omega}, \hat{\mathscr{F}}, \{\hat{\mathscr{F}}_t\}_{t\in[0,T]}, \hat{\mP})$ of $(C_T^\infty, \cB, \{\bar{\cB}_t\}_{t\in[0,T]}, \boldsymbol{\eta})$ such that $\{(\hat{\Omega}, \hat{\mathscr{F}}, \{\hat{\mathscr{F}}_t\}_{t\in[0,T]},
\hat{\mP}), (\boldsymbol{Z}_t=w_t,\hat{W}_t)\}$ is a weak solution of the following stochastic differential equation on $\mR^\infty$: for $u\in\mN$
\be
\dif \boldsymbol{Z}^u_t=<\cT^{-1}(\boldsymbol{Z}_t),\cL_t\varphi_u>\dif t+<\cT^{-1}(\boldsymbol{Z}_t),\varphi_u h^l(t,\cdot)+\partial_i\varphi_u\sigma^{il}_1(t,\cdot)>\dif\hat{W}^l_t, ~ 0\leq t\leq T.
\label{infidime} 
\ee

{\bf Step 2.} We show that Eq.(\ref{zakaieq0}) has a weak solution on $\cM(\mR^n)$.

Put $\hat{\mu}_t:=\cT^{-1}(\boldsymbol{Z}_t)$, and then $\sL_{\hat{\mu}_t}=\Lambda_t$ and (\ref{infidime}) becomes 
\be
<\hat{\mu}_t,\varphi_u>&=&<\hat{\mu}_0,\varphi_u>+\int_0^t<\hat{\mu}_s,\cL_s\varphi_u>\dif s\no\\
&&+\int_0^t<\hat{\mu}_s,\varphi_u h^l(s,\cdot)+\partial_i\varphi_u\sigma^{il}_1(s,\cdot)>\dif\hat{W}^l_s.
\label{apprsequ}
\ee
Thus, we prove that $(\hat{\mu}_t)$ is a weak solution of Eq.(\ref{zakaieq0}).

Now, we specialize the sequence $\{\varphi_u, u\in\mN\}$ as the dense subset of $C_c^\infty(\mR^n)$ (See \cite[Lemma 6.1]{lsz}). And then for any $\varphi\in C_c^\infty(\mR^n)$, there exists a subsequence $\{\varphi_{u_k}, k\in\mN\}$ such that 
$(\varphi_{u_k}, \partial_i\varphi_{u_k}, \partial_{ij}\varphi_{u_k})$ converges unformly to $(\varphi, \partial_i\varphi, \partial_{ij}\varphi)$ as $k\rightarrow\infty$. Accordingly, $<\hat{\mu}_t,\varphi_{u_k}>\rightarrow<\hat{\mu}_t,\varphi>$ as $k\rightarrow\infty$.
Also note that
\ce
&&\hat{\mE}\left|\int_0^t<\hat{\mu}_s,\cL_s\varphi_{u_k}>\dif s-\int_0^t<\hat{\mu}_s,\cL_s\varphi>\dif s\right|\\
&\leq&\hat{\mE}\int_0^t\left|<\hat{\mu}_s,\cL_s(\varphi_{u_k}-\varphi)>\right|\dif s\leq\hat{\mE}\int_0^t\int_{\mR^n}|\cL_s(\varphi_{u_k}-\varphi)(x)|\hat{\mu}_s(\dif x)\dif s\\
&\leq&\hat{\mE}\int_0^t\int_{\mR^n}\bigg[\left|\partial_i(\varphi_{u_k}-\varphi)(x)b_1^i(s,x)\right|+\frac{1}{2}\left|\partial_{ij}(\varphi_{u_k}-\varphi)(x)(\sigma_0\sigma^*_0)^{ij}(s,x)\right|\\
&&\qquad\qquad +\frac{1}{2}\left|\partial_{ij}(\varphi_{u_k}-\varphi)(x)(\sigma_1\sigma^*_1)^{ij}(s,x)\right|\bigg]\hat{\mu}_s(\dif x)\dif s\\
&\leq&C\|\varphi_{u_k}-\varphi\|_{C^2_c}\hat{\mE}\int_0^T\int_{\mR^n}\Big(|b_1(s,x)|+\|\sigma_0\sigma_0^*(s,x)\|+\|\sigma_1(s,x)\|^2\Big)\hat{\mu}_s(\dif x)\dif s\\
&\rightarrow& 0, \qquad k\rightarrow\infty,
\de
where 
$$
\|\varphi\|_{C^2_c}:=\sup\limits_{x\in\mR^n}|\varphi(x)|+\sup\limits_{i}\sup\limits_{x\in\mR^n}|\partial_i\varphi(x)|+\sup\limits_{i,j}\sup\limits_{x\in\mR^n}|\partial_{ij}\varphi(x)|,
$$ 
and
\ce
&&\hat{\mE}\left|\int_0^t<\hat{\mu}_s,\varphi_{u_k} h^l(s,\cdot)+\partial_i\varphi_{u_k}\sigma^{il}_1(s,\cdot)>\dif\hat{W}^l_s-\int_0^t<\hat{\mu}_s,\varphi h^l(s,\cdot)+\partial_i\varphi\sigma^{il}_1(s,\cdot)>\dif\hat{W}^l_s\right|^2\\
&=&\sum\limits_{l=1}^m\hat{\mE}\int_0^t|<\hat{\mu}_s,(\varphi_{u_k}-\varphi) h^l(s,\cdot)+\partial_i(\varphi_{u_k}-\varphi)\sigma^{il}_1(s,\cdot)>|^2\dif s\\
&\leq&\sum\limits_{l=1}^m C\hat{\mE}\int_0^t\int_{\mR^n}\left[|(\varphi_{u_k}-\varphi)(x)h^l(s,x)|^2+|\partial_i(\varphi_{u_k}-\varphi)(x)\sigma^{il}_1(s,x)|^2\right]\hat{\mu}_s(\dif x)\dif s\\
&\leq&C\|\varphi_{u_k}-\varphi\|^2_{C^2_c}\hat{\mE}\int_0^T\int_{\mR^n}\Big(|h(s,x)|^2+\|\sigma_1(s,x)\|^2\Big)\hat{\mu}_s(\dif x)\dif s\\
&\rightarrow& 0, \qquad k\rightarrow\infty.
\de
Thus, replacing $\varphi_u$ in (\ref{apprsequ}) by $\varphi_{u_k}$ and then taking the limit on two sides of (\ref{apprsequ}), we get (\ref{zakaieq2}). The proof is complete.
\end{proof}

Combing Proposition \ref{sdefpe} and Proposition \ref{supeprinpn}, we have the following superposition principle between Eq.(\ref{zakaieq0}) and Eq.(\ref{fpecoup}).

\bt(One superposition principle on $\cM(\mR^n)$)\label{supeprin}

(i)  The existence of weak solutions $\{(\hat{\Omega}, \hat{\mathscr{F}}, \{\hat{\mathscr{F}}_t\}_{t\in[0,T]},\hat{\mP}), (\hat{\mu}_t,\hat{W}_t)\}$ for Eq.(\ref{zakaieq0}) is equivalent to the existence of weak solutions $(\Lambda_t)_{t\in[0,T]}$ for Eq.(\ref{fpecoup}). Moreover, $\Lambda_t=\sL_{\hat{\mu}_t}$ for any $t\in[0,T]$.

(ii) The uniqueness of weak solutions for the Zakai equation (\ref{zakaieq0}) is equivalent to the uniqueness of weak solutions for Eq.(\ref{fpecoup}). 
\et

\subsection{A superposition principle between Eq.(\ref{zakaieq01}) and Eq.(\ref{fpecoup2})}\label{supeprin2}

In the subsection, we state the superposition between Eq.(\ref{zakaieq01}) and Eq.(\ref{fpecoup2}).

\bt(The other superposition principle on $\cM(\mR^n)$)\label{supeprin2}

(i)  The existence of weak solutions $(\hat{\check{\mu}}_t)_{t\in[0,T]}$ for Eq.(\ref{zakaieq01}) is equivalent to the existence of weak solutions $(\check{\Lambda}_t)_{t\in[0,T]}$ for Eq.(\ref{fpecoup2}). Moreover, 
$\check{\Lambda}_t=\sL_{\hat{\check{\mu}}_t}$ for any $t\in[0,T]$.

(ii) The uniqueness of weak solutions for the Zakai equation (\ref{zakaieq01}) is equivalent to the uniqueness of weak solutions for Eq.(\ref{fpecoup2}). 
\et

Since the proof of the above theorem is similar to that in Theorem \ref{supeprin}, we omit it.

Finally, by Remark \ref{aweaksoluzaka2} and Theorem \ref{supeprin2}, we have the following result.

\bc
Suppose that $(i)-(iii)$ hold. Then Eq.(\ref{fpecoup2}) has a unique weak solution.
\ec

\bigskip

\textbf{Acknowledgements:}

The author would like to thank Professor Xicheng Zhang  and Renming Song for valuable discussions.

\end{document}